\documentclass[12pt]{article}
\usepackage{amssymb,amsfonts} 
\usepackage{amsthm}
\usepackage{color}
\usepackage[english]{babel}
 \usepackage{graphicx} 

\newcounter{theorem}
\newenvironment{theorem}[1][\hspace{-1.0ex}]%
 {\par\addvspace{2mm}\indent\refstepcounter{theorem} \textbf{Theorem~\thetheorem\hspace{1.0ex}{\rm#1}.~}\sl}%
 {\par\addvspace{2mm}\rm}
\newcounter{lemma}0
\newenvironment{lemma}[1][\hspace{-1.0ex}]%
 {\par\addvspace{2mm}\indent\refstepcounter{lemma} \textbf{Lemma~\thelemma\hspace{1.0ex}{\rm#1}.~}\sl}%
 {\par\addvspace{2mm}\rm}
\newcounter{corollary}
 {\par\addvspace{2mm}\indent\refstepcounter{corollary} \textbf{Corollary~\thecorollary\hspace{1.0ex}{\rm#1}.~}\sl}%
 {\par\addvspace{2mm}\rm}
\newcounter{note}
 {\par\addvspace{2mm}\indent\refstepcounter{note} \textbf{Remark\hspace{1.0ex}{\rm#1}.~}}%
 {\par\addvspace{2mm}\rm}
\newcounter{example}
 {\par\addvspace{2mm}\indent\refstepcounter{example} \textbf{Example\hspace{1.0ex}{\rm#1}.~}}%
 {\par\addvspace{2mm}\rm}
\newcounter{problem}
 {\par\addvspace{2mm}\indent\refstepcounter{problem} \textbf{Problem~\theproblem\hspace{1.0ex}{\rm#1}.~}}%
 {\par\addvspace{2mm}\rm}

\newcommand\proofr{\par \textbf{Proof.}~}
\newcommand\proofend{\mbox{}\hfill $\blacktriangle$\par}

\def\H#1{H_{#1}}
\def\HH#1{\overline H_{#1}}
\def\HHa#1{\overline H_{#1}^{even}}
\def\HHb#1{\overline H_{#1}^{odd}}

\title{On perfect colorings of the halved $24$-cube%
\thanks{\,The research was partially supported by the RFBR grant 08-01-00673-a}}
\author{Denis S. Krotov}
\begin{document}
\renewcommand\today{}
\maketitle
\renewcommand\abstractname{Abstract}
\begin{abstract}
 A vertex coloring of a graph is said to be perfect with parameters $(a_{ij})_{i,j=1}^k$
 if for every $i,j\in\{1,...,k\}$ every vertex of color $i$ is adjacent with exactly $a_{ij}$
 vertices of color $j$.
  We consider the perfect $2$-colorings of the distance-$2$ graph of the $24$-cube $\{0,1\}^{24}$
  with parameters $((20+c,256-c)(c,276-c))$ (i.e., with eigenvalue $20$).
  We prove that such colorings exist for all $c$ from $1$ to $128$ except
  $1$, $2$, $4$, $5$, $7$, $10$, $13$
  and do not exist for $c=1,2,4,5,7$.
\end{abstract}

\section{Introduction}

We study the vertex $2$-colorings of the vertices of the  distance-$2$ graph of the $24$-cube
that are perfect colorings with the eigenvalue $20$.
These parameters are of interest by the following reasons.

At first,
a known problem is the existence of perfect colorings of the $24$-cube
(i.e., its distance-$1$ graph)
with parameters from the list
$((1,23)\linebreak[2](9,15))$, $((2,22)\linebreak[2](10,14))$, $((3,21)\linebreak[2](11,13))$, $((5,19)\linebreak[2](13,11))$, $((7,17)\linebreak[2](15,9))$
(according to \cite{FDF:PerfColRU,FDF:CorrImmBound,FDF:12cube},
the question of the existence of perfect colorings of the $n$-cube
with fixed parameters is closed for $n<24$).
Colorings with such parameters would correspond to colorings of the distance-$2$ graph
with parameters from the class under study
(see Section \ref{s:12} for the connection between perfect colorings of
the distance-$1$ and distance-$2$ graphs).

At second, an interesting fact is combining two different-nature construction
allows to cover a large specter of $121$ parameter sets.

\vspace{1.5ex}Let $G$ be a simple graph; let $I$ be a finite set, whose elements will be called \emph{colors}.
A coloring $T : V (G) \to I$ is called  \emph{perfect} with parameter matrix
$(s_{ij})_{i,j\in I}$,
iff $T$ is surjective and for every colors $i$ and $j$ every vertex of color $i$ has exactly
$s_{ij}$ color-$j$ neighbors.

By $\H{n}$ we denote the hypercube of dimension $n$, or $n$-cube
(the vertices are the binary words of length $n$;
two words are adjacent iff they differ in exactly one position;
the distance between two words is the number of positions in that they are different).
The the distance-$2$ graph of the hypercube will be denoted by $\HH{n}$
(two words are adjacent iff they differ in exactly two positions);
its degree is $n(n-1)/2$.
This graph has two connected components, which are known as \emph{halved $n$-cubes};
we denote them $\HHa{n}$ and $\HHb{n}$;
from the point of view of perfect colorings it is sufficient to consider only one of them.

Usually, we will consider colorings into two colors, of \emph{$2$-colorings};
the parameter matrix will be written as $((a,b)(c,d))$.
Note that, arranging the colors, we can always set $b\geq c$.
If the graph is regular of degree $s$, then a necessary condition
for the existence of a perfect coloring with parameters $((a,b)(c,d))$ is $a+b=c+d=s$.
So, graph's degree is an eigenvalue of the parameter matrix.
The second  eigenvalue is $a-c=d-b$;
we will refer this value as the \emph{eigenvalue}
of the perfect $2$-coloring and of its parameters.
Often, it is convenient to consider a $2$-coloring as the characteristic
function of some set; in this case $1$ will be considered as the first color; $0$, as the second.

In this paper we study admissible parameters
of perfect $2$-colorings of $\HHa{24}$ with the eigenvalue $20$, i.e., parameters of the form
$((20+c,256-c)(c,276-c))$.

We will prove the following:
\begin{theorem}\label{th:main}
Perfect colorings of $\HHa{24}$ with the parameters $((20+c,256-c)(c,276-c))$
exist for $c=3$, $6$, $8$, $9$, $11$, $12$, and all $c$ from $14$ to $128$.
\end{theorem}
\begin{theorem}\label{th:main2}
There are no perfect colorings of $\HHa{24}$ with the parameters  $((20+c,256-c)(c,276-c))$
if $c$ is $1$, $2$, $4$, $5$, or $7$.
\end{theorem}
The values $c=10$ and $13$ remain under the question; however,
the proved facts allow to observe the existence of gaps in the specter of admissible parameters.
Theorems~\ref{th:main} and~\ref{th:main2} will be proved
in Sections~\ref{s:codes} and~\ref{s:nonexi}.
In Section~\ref{s:12} we will show that a perfect coloring of
$\H{n}$ is a perfect coloring of $\HH{n}$
and establish a relation between the parameters of these colorings.
In Section~\ref{s:unif} observe the possibility to combine two different $2$-colorings
(of an arbitrary graph) with the same eigenvalue
provided the supports of one color are disjoint for two colorings.

\section{A connection between perfect colorings of the graphs $\H{n}$ and $\HH{n}$}\label{s:12}

\begin{lemma}\label{l:d1d2}
A perfect coloring  of $\H{n}$ with matrix $S$ is a perfect coloring  of $\HH{n}$ with matrix
$1/2(S^2-nE)$ (where $E$ is the identity matrix).
\end{lemma}
\proofr
Let us consider a perfect coloring $T$ of $\H{n}$ with parameter matrix $S$.
By the color structure $T(M)$ of some set $M$ of vertices of $\H{n}$ we call the collection
from $|I|$ numbers each of them denoting the number of the vertices of the corresponding color
in $M$.
For an arbitrary vertex $v$ of $\H{n}$
the the color structure of $\{ v \}$ consists of zeros in all the positions except $T(v)$,
where the one is.
Denote by $D_1(v)$ and $D_2(v)$ the set of vertices at the distance $1$ and $2$
from $v$, respectively.
By the definition of a perfect coloring we have $T(V_1(v))=S T(\{v\})$.
Summarizing this formula over the neighborhood $D_1(w)$ of some fixed vertex $w$,
we get
$$
\sum_{v\in D_1(w)}\sum_{u\in D_1(v)}T({u}) = \sum_{v\in D_1(w)}T(V_1(v)) = \sum_{v\in D_1(w)} S T(\{v\})
=  S \sum_{v\in D_1(w)} T(\{v\}) = S^2 T(\{w\}).
$$

On the other hand,
in the sum in the left the index $u$ runs over twice the vertices of $D_2(w)$ and $n$ times,
the vertex $w$. Therefore, this sum also equals $nT({w})+2T(D_2(w))$; so, we deduce
$$
T(D_2(w)) = 1/2 (S^2 T(\{w\}) - n T(\{w\})),
$$
which proves the statement.
\proofend
Note that each of the components of $\HH{n}$ can be colored in less then all colors;
i.e., $\HHa{n}$ (as well as $\HHb{n}$) can be colored in lesser number of colors than $\H{n}$
(see, e.g., Lemma~\ref{l:Gol}).

The following table lists all admissible parameters of perfect $2$-colorings of
$\H{24}$ that correspond to perfect $2$-colorings of $\HH{24}$ with the eigenvalue $\lambda=20$
($20 = 1/2(8^2-24) = 1/2((-8)^2-24)$).
The existence of colorings with the parameters marked by grey
is an open question; for the other parameters, perfect colorings exist
\cite{FDF:PerfColRU,FDF:12cube}.
$$
\small
\def\frm{@{}c@{\,}c@{}}
\def\fra{@{}c@{\ }c@{}}
\begin{array}{@{}r|@{\hspace{0.2ex}}c@{\hspace{0.2ex}}|@{\hspace{0.2ex}}c@{\hspace{0.2ex}}|@{\hspace{0.2ex}}c@{\hspace{0.2ex}}|@{\hspace{0.2ex}}c@{\hspace{0.2ex}}|@{\hspace{0.2ex}}c@{\hspace{0.2ex}}|@{\hspace{0.2ex}}c@{\hspace{0.2ex}}|@{\hspace{0.2ex}}c@{\hspace{0.2ex}}|@{\hspace{0.2ex}}c@{\hspace{0.2ex}}|@{\hspace{0.2ex}}c@{\hspace{0.2ex}}|@{\hspace{0.2ex}}c@{\hspace{0.2ex}}|@{\hspace{0.2ex}}c@{\hspace{0.2ex}}|@{\hspace{0.2ex}}c@{}}
\hline
\begin{array}{@{}r@{}} \H{24} \\ \scriptstyle\lambda{=}{-}8 \end{array}&
&&&
\begin{array}{\fra} 0 & 24 \\ 8 & 16 \end{array}&
\textcolor[gray]{0.6}{\begin{array}{\fra} 1 & 23 \\  9 & 15 \end{array}}&
\textcolor[gray]{0.6}{\begin{array}{\fra} 2 & 22 \\ 10 & 14 \end{array}}&
\textcolor[gray]{0.6}{\begin{array}{\fra} 3 & 21 \\ 11 & 13 \end{array}}&
\begin{array}{\fra} 4 & 20 \\ 12 & 12 \end{array}&
\textcolor[gray]{0.6}{\begin{array}{\fra} 5 & 19 \\ 13 & 11 \end{array}}&
\begin{array}{\fra} 6 & 18 \\ 14 & 10 \end{array}&
\textcolor[gray]{0.6}{\begin{array}{\fra} 7 & 17 \\ 15 &  9 \end{array}}&
\begin{array}{\fra} 8 & 16 \\ 16 & 8 \end{array}
\\ \hline
\begin{array}{@{}r@{}} \H{24} \\ \scriptstyle\lambda{=}8 \end{array}&
\begin{array}{\fra} 9 & 15 \\ 1 & 23 \end{array}&
\begin{array}{\fra} 10 & 14 \\ 2 & 22 \end{array}&
\begin{array}{\fra} 11 & 13 \\ 3 & 21 \end{array}&
\begin{array}{\fra} 12 & 12 \\ 4 & 20 \end{array}&
&
\begin{array}{\fra} 13 & 11 \\ 5 & 19 \end{array}&
&
\begin{array}{\fra} 14 & 10 \\ 6 & 18 \end{array}&
&
\begin{array}{\fra} 15 & 9 \\ 7 & 17 \end{array}&
&
\begin{array}{\fra} 16 & 8 \\ 8 & 16 \end{array} \\ \hline
\begin{array}{@{}r@{}} \HH{24} \\ \scriptstyle\lambda{=}20 \end{array}&
\begin{array}{\frm}  36 & 240 \\  16 & 260 \end{array}&
\begin{array}{\frm}  52 & 224 \\  32 & 244 \end{array}&
\begin{array}{\frm}  68 & 208 \\  48 & 228 \end{array}&
\begin{array}{\frm}  84 & 192 \\  64 & 212 \end{array}&
\begin{array}{\frm}  92 & 184 \\  72 & 204 \end{array}&
\begin{array}{\frm} 100 & 176 \\  80 & 196 \end{array}&
\begin{array}{\frm} 108 & 168 \\  88 & 188 \end{array}&
\begin{array}{\frm} 116 & 160 \\  96 & 180 \end{array}&
\begin{array}{\frm} 124 & 152 \\ 104 & 172 \end{array}&
\begin{array}{\frm} 132 & 144 \\ 112 & 164 \end{array}&
\begin{array}{\frm} 140 & 136 \\ 120 & 156 \end{array}&
\begin{array}{\frm} 148 & 128 \\ 128 & 148 \end{array}\\ \hline
\end{array}
$$

\section{Unifying two perfect colorings with common eigenvalue}
\label{s:unif}
The following lemma, which is straightforward from the definitions,
allow to unify disjoint supports of colors of different colorings with the same eigenvalue.
\begin{lemma}\label{l:union}
Let $C_1$ and $C_2$ be two disjoint subsets of the vertex set $V(G)$
of a simple regular graph $G$; and let $C_1 \cup C_2 \neq V(G)$.
Assume, that the characteristic functions $\chi_{C_1}$ and $\chi_{C_2}$
of $C_1$ and $C_2$ are
perfect colorings of $G$ with the same eigenvalue $\lambda$,
i.e., with the parameters of type
$(\lambda+i,s-\lambda-i;i,s-i)$ и $(\lambda+j,s-\lambda-j;j,s-j)$
where  $s$ is graph's degree.
Then the characteristic function $\chi_{C_1 \cup C_2}$ of the union
is a
perfect colorings of  $G$
with the parameters $(\lambda+i+j,s-\lambda-i-j;i+j,s-i-j)$.
\end{lemma}

\section[Codes and colorings. Proof of Theorem~1]{Codes and colorings. Proof of Theorem~\ref{th:main}}
\label{s:codes}
In this section we construct a class of perfect $2$-colorings of $\HHb{24}$ with the parameters
announced in Theorem~\ref{th:main}.
The support of the first color of a coloring will be constructed
as the union of cosets of one linear code and the neighborhoods of
cosets of the Goley code.

The set of the binary $n$-words (i.e., $V(\H{n})$) will be denoted by $E^n$
and considered as an $n$-dimensional  vector space over the two-element field
with the modulo $2$ calculations.
The \emph{distance} $\rho(\cdot,\cdot)$ between two words is, as usual,
the number of positions in which these words differ
(which coincides with the natural graph metric in $\H{n}$).
Recall that, by definition, an \emph{$(n,M,d)$ code} is a set from $M$ vertices of $\H{n}$
such that the distance between any two different words is not less than $d$.
The \emph{neighborhood} $\Omega(C)$ of some set  $C\subset E^n$
is the set of all the words at the distance $1$ from $C$.

  Let $C_8$ and $C_8'$ be two $(8,16,4)$ codes such that $C_8 \cap C_8' = \{00000000,11111111\}$
  (for definiteness, $C_8$ and $C_8'$ can be defined as containing $00101110$ and $01001110$
  respectively and closed with respect to the addition and  with respect to
   the cyclic permutation of the first seven coordinates).
Define the code
\begin{eqnarray}
\label{eq:Goley}
F=\{(x+y,x+z,x+y+z)|x \in C_8, y,z \in C_8' \}.
\end{eqnarray}
The \emph{distance coloring} of a code $C$ is a function on $E^n$
defined as the distance between the given vertex and $C$.

\begin{lemma}\label{l:Gol}
The distance coloring of the code $F$ is a perfect coloring
 of $\H{24}$ and $\HH{24}$ with matrixes, respectively,
 $$
 \left(
 \begin{array}{ccccc}
   0 & 24 & 0 & 0 & 0 \cr
   1 & 0 & 23 & 0 & 0 \cr
   0 & 2 & 0 & 22 & 0 \cr
   0 & 0 & 3 & 0 & 21 \cr
   0 & 0 & 0 & 24 & 0
 \end{array}
 \right)
\qquad \mbox{and} \qquad
 \left(
 \begin{array}{ccccc}
   0 & 0 & 276 & 0 & 0 \cr
   0 & 23 & 0 & 253 & 0 \cr
   1 & 0 & 44 & 0 & 231 \cr
   0 & 3 & 0 & 273 & 0 \cr
   0 & 0 & 36 & 0 & 240
 \end{array}
 \right).
 $$
 In particular, the generated coloring of $\HHb{24}$
 is a $2$-coloring with parameters $((23,253)(3,273))$.
\end{lemma}
\proofr
It is known \cite[18.7.4]{MWS} that (\ref{eq:Goley}) defines an extended perfect
$(24,2^{12},8)$ code, the Goley code.
This means that the distance from any vertex to $F$ is not more than $4$.
Moreover, the words with even number of ones (from $V(\HHa{24})$) have colors $0$, $2$, $4$;
with odd (from $V(\HHb{24})$), colors $1$, $3$.
So, a color-$4$ vertex is adjacent in $\H{24}$ with $24$ color-$3$ vertices;
a color-$3$ vertex is adjacent with color-$2$ and color-$4$ vertices only,
the number of neighbor color-$2$ vertices being $3$,
because there is exactly one code vertex at the distance $3$ from the given vertex;
the other colors can be checked similarly.
The parameters of the perfect coloring of $\HH{24}$ follow from Lemma~\ref{l:d1d2}.
\proofend

So, $\chi_{\Omega(F)}$ is the first $2$-coloring from the parameter series of Theorem~\ref{th:main}.
Taking disjoint translations of ${\Omega(F)}$ and using Lemma~\ref{l:union},
we would be able to construct $2$-colorings with other parameters.
In order to do it, we need as much as possible
cosets by $F$ at the mutual distance $4$ from each other. Consider the set
$$D=\{(x+y,x+z,x+y+z)|x \in C_8, y,z \in B_8 \},$$
where $B_8$ is the $(8,128,2)$-code containing $00000000$ and, consequently,
including $C_8'$ ($C_8$ and $C_8'$ are defined before the definition of $F$ (\ref{eq:Goley})).

\begin{lemma}\label{l:pho4}
  The set $D$ is a $(24,2^{18},4)$ code.
\end{lemma}
\proofr We first observe the validness of the three simple inequalities
\begin{eqnarray}
\label{eq:uvw}\rho((u,v,w),(u',v',w'))&\geq&\rho(u+v+w,u'+v'+w'),\\
\label{eq:vw}\rho((u,v,w),(u',v',w'))&\geq&\rho(u,u')+\rho(v+w,v'+w'),\\
\label{eq:uw}\rho((u,v,w),(u',v',w'))&\geq&\rho(v,v')+\rho(u+w,u'+w').
\end{eqnarray}
Consider words $r=(x+y,x+z,x+y+z)$ and $r'=(x'+y',x'+z',x'+y'+z')$, where
$x,x'\in C_8$, $y,y',z,z' \in B_8$.
If $x\neq x'$, then, using (\ref{eq:uvw}) and the code distance $4$ of $C_8$, we get
$\rho(r,r') \geq \rho(x,x') \geq 4$. If $x = x'$ and $y \neq y'$, then (\ref{eq:vw})
implies $\rho(r,r') \geq \rho(x+y,x'+y')+\rho(y,y') = 2\rho(y,y') \geq 4$.
The case $x = x'$, $z \neq z'$ is similar.
So, different choices of $x$, $y$, and $z$  lead to different words
with the mutual distance at least $4$ from each other. The number of all such the words is
$|C_8|\cdot|B_8|^2 = 2^{18}$.
\proofend

Since $F$ and $D$ are linear subspaces and, obviously, $F \subset D$,
we can partition $D$ into $64$ cosets by $F$; denote them $F_1,F_2,\ldots,F_{64}$.
Lemma~\ref{l:pho4} implies that the neighborhoods of these cosets are mutually disjoint;
therefore, applying Lemma~\ref{l:union}, we can construct perfect colorings
of $\HHb{24}$ with parameters of type
$((20+3i,256-3i)(3i,276-3i))$, $i=1,\ldots,64$.
In order to cover the larger specter of parameters,
we will need one more code:
$$L=\{(x,y,y+z)|x,y \in B, z \in C_1 \}$$
\begin{lemma}\label{l:8}
  The characteristic function $\chi_L$ is a perfect colorinf of $\HHa{24}$
  with parameters $((28,248)(8,268))$.
\end{lemma}
\proofr
Let as represent $L$ as
$$L=\{(x,w)|x \in B, w \in C_{16} \}$$
where
$C_{16}=\{(y,y+z)|y \in B, z \in C_1 \}$ is a $(16,2^{11},4)$ code.

1) Consider a code vertex $(x,w)$ from $C_{16}$.
The words of $C_{16}$ adjacent with $(x,w)$ in $\HHa{24}$ have the type $(x+e,w)$,
where $e$ is an arbitrary word with exactly two ones. Since the number of such the words $e$ is $28$,
every code vertex is adjacent with exactly $28$ code vertices and, consequently,
 with $248$ non-code vertices.

2) Consider a non-code vertex $(x,w)$ of $\HHa{24}$.
If $x \not\in B$, then the code vertices of $\HHa{24}$
that are adjacent with $(x,w)$ are of type $(x+e,w+e')$ where each of the words $e$ and $e'$
has exactly one one. The word $e$ can be chosen in $8$ ways,
while $e'$, in not more than one way
(otherwise $C_{16}$ contains two words at the distance $2$ from each other).
So, the number of code vertices that are adjacent with $(x,w)$ does not exceed $8$.
If $x \in B$, then $w \not\in C_{16}$, and the code vertices adjacent with $(x,w)$
are of type $(x,w+e'')$ where $e''$ has exactly two ones.
The number of ways to choose $e''$ is not more than $8$
(otherwise $C_{16}$ contains two words at the distance $2$ from each other).
So, every non-code vertex is adjacent with not more than $8$ code ones.
On the other hand, as follows from 1), the number of edges connecting
code and non-code vertices equals
 $2^{18}\cdot 248$ (where $2^{18}$ the number of the code vertices),
 which coincides with $(2^{23}-2^{18})\cdot 8$,
where  $2^{23}-2^{18}$ the number of the non-code vertices ща $\HHa{24}$.
We conclude that every non-code vertex is adjacent with exactly $8$ code vertices and,
consequently, with $268$ non-code ones.
\proofend

Let us consider the set
$$N=\{(x+00000001,y+00000001,z+00000001)|x,y,z \in B)\subset\HHb{24}$$
and partition it into the $8$ cosets $L_1$, \ldots, $L_8$ by $L$.
Since the distance from $D$ to $N$ is $3$, we see that all the sets
 $\Omega(F_1)$, $\Omega(F_2)$, \ldots, $\Omega(F_{64})$, $L_1$, \ldots, $L_8$
 are mutually disjoint and, applying Lemma~\ref{l:union}, get the following:

 \begin{lemma}
   For any $i\in\{0,1, \ldots, 64\}$, $j\in\{0,1, \ldots, 8\}$, $0<i+j<72$,
  the characteristic function of the union of $i$ sets from
  $\Omega(F_1)$, $\Omega(F_2)$, \ldots, $\Omega(F_{64})$
  and $j$ sets from $L_1$, \ldots, $L_8$ is a perfect coloring with parameters
  $((20+3i+8j, 256-3i-8j)( 3i+8j, 276-3i-8j))$.
 \end{lemma}
 Since all the numbers from $1$ to $128$ except $1$, $2$, $4$, $5$, $7$, $10$, $13$
 can be represented as $3i+8j$, Theorem~\ref{th:main} is proved.

 \section[Proof of Theorem~2. The nonexistence]{Proof of Theorem~\ref{th:main2}. The nonexistence}
\label{s:nonexi}
In this section we prove the nonexistence
of perfect colorings with parameters
$((20+c,256-c)(c,276-c))$ for $c=1,2,4,5,7$.
A set $S \subset V(\HHa{24})$ is called a \emph{sphere} iff it consists of all the $24$ vertices
at the distance $1$ from some fixed vertex in $V(\HHb{24})$.

In the following proof the notion of type $V_{24}$, where $V\subset \{1,\ldots,24\}$,
means the binary length-$24$ word with the nonzero-position set $V$.

\begin{lemma}\label{l:12345}
Assume that the characteristic function $\chi_C$
of a set $C \subset V(\HHa{24})$
is a perfect coloring of $\HHa{24}$ with parameters $((20+c,256-c)(c,276-c))$.
If $c\leq 7$, then $C$ is the union of spheres.
\end{lemma}
\proofr
We consider only the two cases of $((25,251)(5,271))$ and $((27,249)(7,269))$,
because the other cases are proved similarly.

In the first case, $c=5$.
Let us take an arbitrary $v$ from $C$ and show that it belongs to a sphere included to $C$.
W.l.o.g. we can assume
$v=000000000000000000000000$. A pair  $\{i,j\}$ of coordinates from $1$ to $24$ is called
\emph{code} iff $\{i,j\}_{24}\in C$.
As follows from the parameter matrix, there are exactly $25$ code pairs.
Moreover, $c=5$ also implies that

(*) a non-code pair intersects with at most four code pairs
(the fifth neighbor will be $v$).

Consider the cases:

1) If some $i$th coordinate belongs to $23$ code pairs, then the corresponding words
together with $v$ constitute a sphere, which proves the statement for this case.

2) If some $i$th coordinate belongs to less than $23$ and more than $4$ code pairs,
then
$i$ belongs to some non-code pair, contradicting (*).

3) If some $i$th coordinate belongs to exactly $4$ code pairs, then there is a
code pair $\{j,k\}$ disjoint with all of them.
Then the non-code pair $\{i,j\}$ contradicts (*).

4) Assume that there is no a coordinate that belongs to more than $3$ code pairs.
Since the number of code pairs is greater than $24$, there is a coordinate $i$
that belongs to some $3$ code pairs
$\{i,j_1\}$, $\{i,j_2\}$, $\{i,j_3\}$. It is easy to count that
among the remaining $20$ coordinates there is $j$ that belongs to two code pairs.
Then the non-code pair $\{i,j\}$ contradicts (*).

For the case $c=5$, the claim of the lemma is proved.

Let us consider the case $((27,249)(7,269))$, i.e., $c=7$.
In general, the idea of the proof is similar to the previous considered case,
but now we will consider not one but two neighbor vertices of first color
and show that for any coloring of the neighborhood of
these two vertices the condition $a=27$ contradicts to $c=7$
(assuming that $C$ is not the union of spheres).

By $C'$ we denote the union of the spheres included in $C$;
by $C''$, the set $C \setminus C'$.
We have to prove that $C''$ is empty. Suppose the contrary.
Let $k$ be the maximum cardinality of the intersection of $C$ and
a sphere that contains at least one vertex from $C''$.
Since such a sphere must contain a vertex not from $C$, we have
\begin{equation}\label{eq:k7}
k\leq c=7.
\end{equation}
W.l.o.g. we can assume that
 ${\bar 0}=000000000000000000000000\in C''$ and $\{1,2\}_{24},$ $\{1,3\}_{24},$ \ldots,
 $\{1,k\}_{24}$ $\in C$.

Let us consider three matrices whose rows are words of $C$
from the neighborhoods of ${\bar 0}$ and $\{1,2\}_{24}$.
The rows of the matrix $A_1$ are the weight-$2$ vectors different
from $\{1,2\}_{24}$ and non-adjacent with $\{1,2\}_{24}$.
The weight-$2$ (weight-$4$) vectors adjacent with $\{1,2\}_{24}$
form the matrix $A_2$ (respectively, $A_3$).
The vectors ${\bar 0}$ and $\{1,2\}_{24}$ themselves are not included to one of the matrices.
The matrices $A_1$ and $A_2$ (as well as $A_2$ and $A_3$) have, summarily,
$a-1=26$ rows, which with $\{1,2\}_{24}$ (respectively, ${\bar 0}$) form the
intersection of $C$ and the neighborhood of ${\bar 0}$ (respectively, $\{1,2\}_{24}$).
Denoting the height of $A_i$ by $h_i$, we have
$$h_1+h_2= h_2 + h_3 = 26.$$
Moreover, from the definition of $k$ we deduce $$h_2\leq 2k-4,$$ because the rows of $A_2$
have the form $\{1,j\}_{24}$ (not more than $k-2$ rows, because together with ${\bar 0}$ and $\{1,2\}_{24}$
they belong to the sphere centered in $\{1\}_{24}$) or $\{2,j\}_{24}$ (similarly).

Since every row of $A_1$ and $A_2$ contains exactly two ones,
the total number of ones in $A_1$ and $A_2$ is $52$.
Thus, there is a column that contains (summarily in the two matrices)
at least three ones, which means by the definition of $k$ that
\begin{equation}\label{eq:k4}
k\geq 4.
\end{equation}
Every row of $A_3$ contains four ones. In summary, the number of ones in the three matrices
equals
\begin{equation}\label{eq:111}
2h_1+2h_2+4h_3 = 2\cdot 26 + 2(26-h_2)+2h_3 \geq 104 - 2(2k-4)+2h_3 = 112-4k+2h_3
\end{equation}
($2h_3$ corresponds to the ones in the first two columns of $A_3$
and will be canceled in the following estimation).

Let us estimate this value from the other side.
We call the columns with the numbers at most (more than)
$k$ the left (respectively, right) part of the matrix.

(a) \emph{The left part of $A_1$ contains at most $(k-2)(k-1)$ ones.} Indeed,
the first two columns of $A_1$ are zero; and, by the definition of $k$, any other column
contains at most $k-1$ ones.

(b) \emph{The left part of $A_2$ contains at most $4(k-2)$ ones.} Indeed,
as noted above, $A_2$ has at most
$2(k-2)$ rows, every row containing exactly two ones.

(c) \emph{The right parts of the three matrices contain summarily at most $(24-k)(7-k)$ ones.}
Indeed, if for some $j>k$ the $j$th columns of $A_1$, $A_2$, and $A_3$ contain
more than
$7-k$ ones, then the word $\{1,j\}_{24}\not\in C$ has more than $7$ neighbors from $C$
(the corresponding rows plus
${\bar 0}$, $\{1,2\}_{24}$, $\{1,3\}_{24}$, \ldots, $\{1,k\}_{24}$),
which contradicts to the parameter $c=7$.

Let us consider separately two subcases:

I. $\{1,2\}_{24}\in C''$. In this case we additionally have the following:

(d) \emph{The left part of $A_3$ contains at most $2h_3+(k-2)(k-1)$ ones.}
(Similarly to (a), but the first two columns consist of ones.)

In summary, estimating the total number of ones in $A_1$, $A_2$, and $A_3$ from (a)-(d)
and taking into account $(\ref{eq:111})$, we have
$$112-4k+2h_3 \leq (k-2)(k-1)+4(k-2)+2h_3+(k-2)(k-1)+(24-k)(7-k).$$
I.e.,
$$ 3k^2-29k+52 \geq 0,$$
which does not hold for the values $k=4,5,6,7$ satisfying (\ref{eq:k7}) and (\ref{eq:k4}).
This contradiction proves the statement for the subcase I.

II. $\{1,2\}_{24}\in C'$. Then, $C$ includes a sphere centered in $\{1,2,j\}_{24}$ for some $j$.
Consequently, $A_3$ contains all rows of type $\{1,2,j,i\}_{24}$,
$i\in \{3,\ldots,24\}\setminus\{j\}$; thus, every column of $A_3$ contains a one.
Taking into account (c), we get the following:

(e) \emph{The right parts of $A_1$ and $A_2$ contain at most $(24-k)(6-k)$ ones summarily.}

As noted above, the total number of ones in $A_1$ and $A_2$ is $52$.
On the other hand, as follows from (a), (b), and (e), it is not greater than
$$ (k-2)(k-1) + 4(k-2) + (24-k)(6-k) = 2k^2 - 29k +138. $$
We deduce that $2k^2 - 29k + 86 \geq 0$, which is not true if $5\leq k \leq 7$,
but holds for $k=4$.

Let us consider this remaining subcase.
We have: $\{1,2\}_{24}$, $\{1,3\}_{24}$, $\{1,4\}_{24}$ $\in C$;
moreover, $\{1,2\}_{24}$ belongs to a sphere included in $C$.
The center of the sphere has the form $\{1,2,j\}_{24}$
(it cannot contain only one one, because $\bar 0 \in C''$).
Then, $\{1,j\}_{24}\in C$; so, $j$ is $3$ or $4$.
Assume w.l.o.g. that $j=3$.
We claim that

(**) $\{1,4\}_{24}\in C''$. Suppose, by contradiction, that
$\{1,4\}_{24}$ belongs to a sphere included in $C$.
Similarly to the arguments above, the center of the sphere
must have the form $\{1,j,4\}_{24}$ where $j$ is $2$ or $3$.
But the spheres with centers $\{1,2,3\}_{24}$ and $\{1,j,4\}_{24}$ have nonempty intersection.
A vertex from the intersection belongs to $C$
and has at least $45$ neighbors from $C$
($1+45$ is the cardinality of the union of the two spheres),
which contradicts to $a=27$. The claim (**) is proved.

So, interchanging the second and the fourth coordinates leads to the case I.
\proofend

\begin{lemma}\label{l:::3}
Let $C \subset V(\HHa{24})$ be the union of spheres; and let the characteristic function $\chi_C$
be a perfect coloring of $\HHa{24}$ with $((20+c,256-c)(c,276-c))$.
Then either $c$ is divisible $3$, or $c \geq 25$.
\end{lemma}
\proofr
If $C$ includes two intersecting spheres, then their union contains $46$  vertices,
and a common vertex is adjacent with the other $45$ vertices of the union.
Thus, $20+c \geq 45$, i.e., $c \geq 25$.

If, otherwise, $C$ consists of disjoint spheres,
then every vertex from $V(\HHa{24})\setminus C$
is adjacent with exactly three vertices from each neighbor sphere; thus, $c\equiv 0 \bmod 3$.
\proofend
Since the values $c=1,2,4,5,7$ contradict Lemmas~\ref{l:12345} and~\ref{l:::3},
Theorem~\ref{th:main2} is proved.
Moreover, we can conclude that a perfect coloring with parameters $((23,253)(3,273))$
is unique up to graph automorphisms.
\section[Conclusion]{Conclusion}

For the conclusion, we list all values of $c$ from $1$ to $128$ in a table.
The sign "$-$"\ means the nonexistence of perfect colorings with parameters
$(20+c,256-c;c,276-c)$ in $\HH{24}$;
"$+$"{}, the
existence, "$?$"\ the question is open.
A circle means that a coloring with these parameters can be constructed as a union
of cosets of the neighborhood of the Goley code;
a box,
as a union
of cosets of the linear code $L$.

$$\scriptstyle
\def\cpl{\makebox[0mm][c]{$\bigcirc$}\makebox[0mm][c]{$+$}}
\def\qpl{\framebox[4.5mm]{$+$}}
\def\qcpl{\framebox[4.5mm]{\makebox[0mm][c]{$\bigcirc$}\makebox[0mm][c]{$+$}}}
\begin{array}{|@{\,}c@{\,}|@{\,}c@{\,}|@{\,}c@{\,}|@{\,}c@{\,}|@{\,}c@{\,}|@{\,}c@{\,}|@{\,}c@{\,}|@{\,}c@{\,}|@{\,}c@{\,}|@{\,}c@{\,}|@{\,}c@{\,}|@{\,}c@{\,}|@{\,}c@{\,}|@{\,}c@{\,}|@{\,}c@{\,}|@{\,}c@{\,}|@{\,}c@{\,}|@{\,}c@{\,}|@{\,}c@{\,}|@{\,}c@{\,}|@{\,}c@{\,}|@{\,}c@{\,}|@{\,}c@{\,}|@{\,}c@{\,}|@{\,}c@{\,}|@{\,}c@{\,}|@{\,}c@{\,}|@{\,}c@{\,}|@{\,}c@{\,}|@{\,}c@{\,}|@{\,}c@{\,}|@{\,}c@{\,}|}
\hline
1 & 2 & 3 & 4 & 5 & 6 & 7 & 8 & 9 & 10 & 11 & 12 & 13 & 14 & 15 & 16 \\
-&-&\cpl&-&-&\cpl&-&\qpl&\cpl&?&+&\cpl&?&+&\cpl&\qpl\\ \hline
 17 & 18 & 19 & 20 & 21 & 22 & 23 & 24 & 25 & 26 & 27 & 28 & 29 & 30 & 31 & 32 \\
+&\cpl&+&+&\cpl&+&+&\qcpl&+&+&\cpl&+&+&\cpl&+&\qpl\\ \hline
33 & 34 & 35 & 36 & 37 & 38 & 39 & 40 & 41 & 42 & 43 & 44 & 45 & 46 & 47 & 48 \\
\cpl&+&+&\cpl&+&+&\cpl&\qpl&+&\cpl&+&+&\cpl&+&+&\qcpl\\ \hline
 49 & 50 & 51 & 52 & 53 & 54 & 55 & 56 & 57 & 58 & 59 & 60 & 61 & 62 & 63 & 64 \\
+&+&\cpl&+&+&\cpl&+&\qpl&\cpl&+&+&\cpl&+&+&\cpl&\qcpl\\ \hline
65  &   66  &   67  &   68  &   69  &   70  &   71  &   72  &   73  &   74  &   75  &   76  &   77  &   78  &   79  &   80  \\
+&\cpl&\cpl&+&\cpl&\cpl&+&\qcpl&\cpl&+&\cpl&\cpl&+&\cpl&\cpl&\qpl\\ \hline
81  &   82  &   83  &   84  &   85  &   86  &   87  &   88  &   89  &   90  &   91  &   92  &   93  &   94  &   95  &   96  \\
\cpl&\cpl&+&\cpl&\cpl&+&\cpl&\qcpl&+&\cpl&\cpl&+&\cpl&\cpl&+&\qcpl\\ \hline
97  &   98  &   99  &   100 &   101 &   102 &   103 &   104 &   105 &   106 &   107 &   108 &   109 &   110 &   111 &   112 \\
\cpl&+&\cpl&\cpl&+&\cpl&\cpl&\qpl&\cpl&\cpl&+&\cpl&\cpl&+&\cpl&\qcpl\\ \hline
113 &   114 &   115 &   116 &   117 &   118 &   119 &   120 &   121 &   122 &   123 &   124 &   125 &   126 &   127 &   128 \\
+&\cpl&\cpl&+&\cpl&\cpl&+&\qcpl&\cpl&+&\cpl&\cpl&+&\cpl&\cpl&\qpl\\ \hline
\end{array}
$$

\providecommand\href[2]{#2} \providecommand\url[1]{\href{#1}{#1}}

\end{document}


\renewcommand\today{}
\maketitle
\begin{abstract}
 Раскраска вершин графа называется совершенной с параметрами $(a_{ij})_{i,j=1}^k$,
 если для всех $i,j\in\{1,...,k\}$ каждая вершина цвета $i$ смежна ровно с $a_{ij}$
 вершинами цвета $j$.
  Рассматриваются совершенные раскраски в два цвета графа расстояний $2$ гиперкуба
  $\{0,1\}^{24}$
  с параметрами $((20+c,256-c)(c,276-c))$ (т. е., с собственным значением $20$).
  Мы доказываем, что такие раскраски существуют при
  всех $c$ от $1$ до $128$ кроме $1$, $2$, $4$, $5$, $7$, $10$, $13$
  и не существуют при $c=1,2,4,5,7$.
\end{abstract}

\section[Introduction]{Введение}

Мы изучаем  раскраски в два цвета вершин графа расстояний два $24$-мерного гиперкуба,
являющиеся совершенными раскрасками с собственным числом $20$.
Интерес именно к этим параметрам обусловлен следующими причинами.

Во-первых,
известен вопрос существования совершенных раскрасок гиперкуба (то есть его графа расстояний $1$)
размерности $24$ с параметрами
$((1,23)(9,15))$, $((2,22)\linebreak[2](10,14))$, $((3,21)(11,13))$, $((5,19)(13,11))$, $((7,17)(15,9))$
(согласно \cite{FDF:PerfColRU,FDF:CorrImmBound,FDF:12cube},
вопрос существования совершенных раскрасок $n$-куба с фиксированными параметрами полностью исследован для $n<24$).
Раскраски с этими параметрами соответствовали бы раскраскам графа расстояний $2$ с параметрами
из исследуемого класса (подробнее о связи параметров раскрасок графов расстояний $1$ и $2$ см. раздел \ref{s:12}).

Во-вторых, интересен факт, что, комбинируя две совершенно разнородные конструкции,
удаётся покрыть большой спектр из $121$ набора параметров.

\vspace{1.5ex}Пусть $G$ -- простой граф, $I$ -- конечное множество, элементы которого будем называть \emph{цветами}.
Отображение $T : V (G) \to I$ называется \emph{совершенной раскраской} с матрицей параметров
$(s_{ij})_{i,j\in I}$,
если оно сюръективно и для каждых $i$, $j$ у каждой вершины цвета $i$ число соседей
цвета $j$ равно $s_{ij}$.

Через $\H{n}$ обозначим гиперкуб размерности $n$
(вершины \emph{гиперкуба}, или \emph{$n$-куба}, -- двоичные слова длины $n$;
два слова смежны если и только если они различаются ровно в одной позиции;
расстояние между двумя словами есть число позиций, в которых слова различны).
Граф расстояний $2$ гиперкуба будем обозначать через $\HH{n}$
(два слова смежны если и только если они различаются ровно в двух позициях),
его степень равна $n(n-1)/2$.
Поскольку этот граф состоит из двух компонент связности, известные как
\emph{половинчатые} (halved) \emph{$n$-кубы}
обозначим их $\HHa{n}$ и $\HHb{n}$,
с точки зрения совершенных раскрасок достаточно рассматривать одну из компонент.

Как правило, мы будем рассматривать раскраски в два цвета, или \emph{$2$-раскраски};
матрицу параметров будем записывать в виде $((a,b)(c,d))$.
Заметим, что, нужным образом упорядочив цвета, всегда можно добиться $b\geq c$.
Если граф регулярный степени $s$, необходимыми условием существования совершенной раскраски
с параметрами $((a,b)(c,d))$ является $a+b=c+d=s$.
Таким образом, степень графа
является собственным числом матрицы параметров.
Вторым собственным числом матрицы является $a-c=d-b$,
это значение будем считать \emph{собственным числом}
совершенной раскраски и её параметров.
Часто удобно считать $2$-раскраску характеристической функцией некоторого множества,
в этом случае первым цветом договоримся считать $1$, вторым -- $0$.

В данной работе исследуются возможные параметры совершенных раскрасок в два цвета
графа $\HHa{24}$ с собственным числом $20$, т.\,е. параметры вида $((20+c,256-c)(c,276-c))$.
Результатом являются следующие два утверждения.
\begin{theorem}\label{th:main}
Совершенные раскраски графа $\HHa{24}$ с параметрами вида $((20+c,256-c)(c,276-c))$
существуют при $c=3$, $6$, $8$, $9$, $11$, $12$ и всех $c$ от $14$ до $128$.
\end{theorem}
\begin{theorem}\label{th:main2}
Не существует совершенных раскрасок графа $\HHa{24}$ с параметрами вида $((20+c,256-c)(c,276-c))$
при $c$ равном $1$, $2$, $4$, $5$ или $7$.
\end{theorem}
Значения $10$ и $13$ остаются под вопросом, однако, доказанное уже позволяет
отметить разрывы в спектре допустимых значений.
Теоремы~\ref{th:main} и~\ref{th:main2} будут доказаны в разделах~\ref{s:codes} и~\ref{s:nonexi}.
В разделе~\ref{s:12} показано, что совершенная раскраска графа
$\H{n}$ является совершенной раскраской графа $\HH{n}$
и установлена связь параметров этих раскрасок.
В разделе~\ref{s:unif} указано на возможность комбинирования $2$-раскрасок
(произвольного графа) с одинаковым собственным числом,
при условии непересекаемости носителей одного из цветов.

\section[A connection between perfect colorings of the graphs $\H{n}$ and $\HH{n}$]%
{Связь совершенных раскрасок графов $\H{n}$ и $\HH{n}$}\label{s:12}

\begin{lemma}\label{l:d1d2}
Совершенная раскраска с матрицей $S$ графа $\H{n}$ является совершенной раскраской
с матрицей $1/2(S^2-nE)$ графа $\HH{n}$ на том же множестве вершин ($E$ -- единичная матрица).
\end{lemma}
\proofr
Рассмотрим совершенную раскраску $T$ графа $\H{n}$ с матрицей параметров $S$.
Цветовым составом $T(M)$
множества вершин $M$ графа $\H{n}$ назовём набор из $|I|$ чисел
каждое из которых обозначает число вершин соответствующего цвета в множестве $M$.
Для произвольной вершины $v$ графа $\H{n}$
цветовой состав множества $\{ v \}$ состоит из нулей во всех позициях, кроме $T(v)$,
где стоит единица.
Обозначим через $D_1(v)$ и $D_2(v)$ множество вершин на расстоянии $1$ и $2$
от $v$, соответственно.
По определению совершенной раскраски имеем $T(V_1(v))=S T(\{v\})$.
Посуммировав эту формулу по окрестности $D_1(w)$ некоторой фиксированной вершины $w$,
получим
$$
\sum_{v\in D_1(w)}\sum_{u\in D_1(v)}T({u}) = \sum_{v\in D_1(w)}T(V_1(v)) = \sum_{v\in D_1(w)} S T(\{v\})
=  S \sum_{v\in D_1(w)} T(\{v\}) = S^2 T(\{w\}).
$$

С другой стороны, в сумме слева индекс $u$ дважды пробегает вершины $D_2(w)$ и $n$ раз -- саму
вершину $w$. Таким образом, эта сумма равна также $nT({w})+2T(D_2(w))$, откуда
$$
T(D_2(w)) = 1/2 (S^2 T(\{w\}) - n T(\{w\})),
$$
что и требовалось доказать.
\proofend
Отметим, что может оказаться, что в каждой компоненте связности графа $\HH{n}$
будут встречаться не все цвета, то есть граф $\HHa{n}$ (как и $\HHb{n}$)
будет раскрашен в меньшее число цветов (см. например Лемму~\ref{l:Gol}).

В следующей таблице приведены все возможные параметры совершенных раскрасок в два цвета
графа $\H{24}$, которые соответствуют совершенным раскраскам $\HH{24}$ с собственным числом $\lambda=20$
($20 = \frac{8^2-24}2 = \frac{(-8)^2-24}2$).
Про существование раскрасок с параметрами, затенёнными серым цветом,
в настоящее время ничего не известно,
раскраски с остальными параметрами существуют \cite{FDF:PerfColRU}, \cite{FDF:12cube}.
$$
\small
\def\frm{@{}c@{\,}c@{}}
\def\fra{@{}c@{\ }c@{}}
\begin{array}{@{}r|@{\hspace{0.2ex}}c@{\hspace{0.2ex}}|@{\hspace{0.2ex}}c@{\hspace{0.2ex}}|@{\hspace{0.2ex}}c@{\hspace{0.2ex}}|@{\hspace{0.2ex}}c@{\hspace{0.2ex}}|@{\hspace{0.2ex}}c@{\hspace{0.2ex}}|@{\hspace{0.2ex}}c@{\hspace{0.2ex}}|@{\hspace{0.2ex}}c@{\hspace{0.2ex}}|@{\hspace{0.2ex}}c@{\hspace{0.2ex}}|@{\hspace{0.2ex}}c@{\hspace{0.2ex}}|@{\hspace{0.2ex}}c@{\hspace{0.2ex}}|@{\hspace{0.2ex}}c@{\hspace{0.2ex}}|@{\hspace{0.2ex}}c@{}}
\hline
\begin{array}{@{}r@{}} \H{24} \\ \scriptstyle\lambda{=}{-}8 \end{array}&
&&&
\begin{array}{\fra} 0 & 24 \\ 8 & 16 \end{array}&
\textcolor[gray]{0.6}{\begin{array}{\fra} 1 & 23 \\  9 & 15 \end{array}}&
\textcolor[gray]{0.6}{\begin{array}{\fra} 2 & 22 \\ 10 & 14 \end{array}}&
\textcolor[gray]{0.6}{\begin{array}{\fra} 3 & 21 \\ 11 & 13 \end{array}}&
\begin{array}{\fra} 4 & 20 \\ 12 & 12 \end{array}&
\textcolor[gray]{0.6}{\begin{array}{\fra} 5 & 19 \\ 13 & 11 \end{array}}&
\begin{array}{\fra} 6 & 18 \\ 14 & 10 \end{array}&
\textcolor[gray]{0.6}{\begin{array}{\fra} 7 & 17 \\ 15 &  9 \end{array}}&
\begin{array}{\fra} 8 & 16 \\ 16 & 8 \end{array}
\\ \hline
\begin{array}{@{}r@{}} \H{24} \\ \scriptstyle\lambda{=}8 \end{array}&
\begin{array}{\fra} 9 & 15 \\ 1 & 23 \end{array}&
\begin{array}{\fra} 10 & 14 \\ 2 & 22 \end{array}&
\begin{array}{\fra} 11 & 13 \\ 3 & 21 \end{array}&
\begin{array}{\fra} 12 & 12 \\ 4 & 20 \end{array}&
&
\begin{array}{\fra} 13 & 11 \\ 5 & 19 \end{array}&
&
\begin{array}{\fra} 14 & 10 \\ 6 & 18 \end{array}&
&
\begin{array}{\fra} 15 & 9 \\ 7 & 17 \end{array}&
&
\begin{array}{\fra} 16 & 8 \\ 8 & 16 \end{array} \\ \hline
\begin{array}{@{}r@{}} \HH{24} \\ \scriptstyle\lambda{=}20 \end{array}&
\begin{array}{\frm}  36 & 240 \\  16 & 260 \end{array}&
\begin{array}{\frm}  52 & 224 \\  32 & 244 \end{array}&
\begin{array}{\frm}  68 & 208 \\  48 & 228 \end{array}&
\begin{array}{\frm}  84 & 192 \\  64 & 212 \end{array}&
\begin{array}{\frm}  92 & 184 \\  72 & 204 \end{array}&
\begin{array}{\frm} 100 & 176 \\  80 & 196 \end{array}&
\begin{array}{\frm} 108 & 168 \\  88 & 188 \end{array}&
\begin{array}{\frm} 116 & 160 \\  96 & 180 \end{array}&
\begin{array}{\frm} 124 & 152 \\ 104 & 172 \end{array}&
\begin{array}{\frm} 132 & 144 \\ 112 & 164 \end{array}&
\begin{array}{\frm} 140 & 136 \\ 120 & 156 \end{array}&
\begin{array}{\frm} 148 & 128 \\ 128 & 148 \end{array}\\ \hline
\end{array}
$$

\section[Unifying two perfect colorings with common eigenvalue]{Объединение цветов двух раскрасок с общим собственным числом}
\label{s:unif}
Следующая лемма, вытекающая непосредственно из определений, позволяет объединять
непересекающиеся носители цветов различных раскрасок с одинаковым собственным числом.
\begin{lemma}\label{l:union}
Пусть $C_1$ и $C_2$ -- два непересекающихся подмножества множества вершин $V(G)$
простого регулярного графа $G$, причём $C_1 \cup C_2 \neq V(G)$.
И предположим, что характеристические функции $\chi_{C_1}$ и $\chi_{C_2}$
множеств $C_1$ и $C_2$ являются
совершенными раскрасками $G$ с одинаковым собственным числом $\lambda$,
то есть с параметрами вида
$(\lambda+i,s-\lambda-i;i,s-i)$ и $(\lambda+j,s-\lambda-j;j,s-j)$
где $s$ -- степень графа.
Тогда характеристическая функция $\chi_{C_1 \cup C_2}$ объединения
есть совершенная раскраска $G$
с параметрами $(\lambda+i+j,s-\lambda-i-j;i+j,s-i-j)$.
\end{lemma}

%

\section[Codes and colorings. Proof of Theorem~1]{Коды и раскраски. Доказательство
теоремы~\ref{th:main}}
\label{s:codes}
В этом разделе мы построим класс совершенных раскрасок графа $\HHb{24}$ с параметрами,
анонсированными в теореме~\ref{th:main}.
Носитель первого цвета будет строиться как объединение смежных классов
одного линейного кода и окрестностей смежных классов кода Голея.

Множество двоичных слов длины $n$ (т.\,е. $V(\H{n})$)
будем обозначать через $E^n$ и рассматривать
как $n$-мерное векторное пространство над полем из двух элементов,
с вычислениями по модулю $2$.
Под \emph{расстоянием} $\rho(\cdot,\cdot)$ между двумя словами будем, как обычно, подразумевать
число позиций, в которых эти два слова различны (что совпадает с естественной метрикой
в графе $\H{n}$).
Напомним, что \emph{$(n,M,d)$-кодом} называется множество из $M$ вершин $\H{n}$ такое, что
расстояние между любыми двумя различными вершинами не меньше $d$.
\emph{Окрестностью} $\Omega(C)$ некоторого множества $C$
будем считать множество слов на расстоянии $1$ от $C$.

  Пусть $C_8$ и $C_8'$ -- два $(8,16,4)$-кода таких, что $C_8 \cap C_8' = \{00000000,11111111\}$
  (для определённости, $C_8$ и $C_8'$ можно определить как содержащие $00101110$ и $01001110$ соответственно
  и замкнутые относительно сложения и циклической перестановки первых семи координат).
Определим код
\begin{eqnarray}
\label{eq:Goley}
F=\{(x+y,x+z,x+y+z)|x \in C_8, y,z \in C_8' \}.
\end{eqnarray}
\emph{Дистанционной раскраской} кода $C$ назовём функцию,
сопоставляющую каждой вершине из $\H{n}$
расстояние от неё до $C$.

\begin{lemma}\label{l:Gol}
 Дистанционная раскраска кода $F$ есть совершенная раскраска
 графов $\H{24}$ и $\HH{24}$ с матрицами, соответственно,
 $$
 \left(
 \begin{array}{ccccc}
   0 & 24 & 0 & 0 & 0 \cr
   1 & 0 & 23 & 0 & 0 \cr
   0 & 2 & 0 & 22 & 0 \cr
   0 & 0 & 3 & 0 & 21 \cr
   0 & 0 & 0 & 24 & 0
 \end{array}
 \right)
\qquad \mbox{и} \qquad
 \left(
 \begin{array}{ccccc}
   0 & 0 & 276 & 0 & 0 \cr
   0 & 23 & 0 & 253 & 0 \cr
   1 & 0 & 44 & 0 & 231 \cr
   0 & 3 & 0 & 273 & 0 \cr
   0 & 0 & 36 & 0 & 240
 \end{array}
 \right).
 $$
 В частности, порождённая раскраска графа $\HHb{24}$
 есть $2$-раскраска с параметрами $((23,253)(3,273))$.
\end{lemma}
\proofr
Известно \cite[18.7.4]{MWS}, что (\ref{eq:Goley}) задаёт расширенный совершенный
$(24,2^{12},8)$-код, известный как код Голея.
Это означает, что расстояние от любой вершины до кода не больше четырёх.
Кроме того, вершины с чётным числом единиц (из $V(\HHa{24})$) имеют цвета $0$, $2$, $4$,
а с нечётным (из $V(\HHb{24})$)  -- цвета $1$, $3$.
Таким образом, вершина цвета $4$ смежна в $\H{24}$ с $24$ вершинами цвета $3$;
вершина цвета $3$ смежна только с вершинами цветов $2$ и $4$,
причём смежных вершин цвета $2$
ровно $3$, поскольку на расстоянии $3$ от данной находится только одна кодовая вершина;
аналогично проверяются другие цвета.
Параметры раскраски графа $\HH{24}$ следуют из леммы~\ref{l:d1d2}.
\proofend

Таким образом, $\chi_{\Omega(F)}$ есть первая совершенная $2$-раскраска из серии,
описанной теоремой~\ref{th:main}. Взяв непересекающиеся сдвиги ${\Omega(F)}$ и
воспользовавшись леммой~\ref{l:union}, мы можем строить раскраски с другими параметрами.
Для этого нам нужно набрать как можно больше смежных классов по $F$ на расстоянии не менее $4$
друг от друга. Рассмотрим множество
$$D=\{(x+y,x+z,x+y+z)|x \in C_8, y,z \in B_8 \},$$
где $B$ есть $(8,128,2)$-код, содержащий $00000000$ и, следовательно, включающий
$C_8'$ ($C_8$ и $C_8'$ вводились перед определением кода $F$ (\ref{eq:Goley})).

\begin{lemma}\label{l:pho4}
  Множество $D$ есть $(24,2^{18},4)$-код.
\end{lemma}
\proofr Для начала отметим три простых неравенства
\begin{eqnarray}
\label{eq:uvw}\rho((u,v,w),(u',v',w'))&\geq&\rho(u+v+w,u'+v'+w'),\\
\label{eq:vw}\rho((u,v,w),(u',v',w'))&\geq&\rho(u,u')+\rho(v+w,v'+w'),\\
\label{eq:uw}\rho((u,v,w),(u',v',w'))&\geq&\rho(v,v')+\rho(u+w,u'+w').
\end{eqnarray}
Рассмотрим слова $r=(x+y,x+z,x+y+z)$ и $r'=(x'+y',x'+z',x'+y'+z')$, где
$x,x'\in C_8$, $y,y',z,z' \in B_8$.
Если $x\neq x'$, то, воспользовавшись (\ref{eq:uvw}) и расстоянием $4$ кода $C_8$, получаем
$\rho(r,r') \geq \rho(x,x') \geq 4$. Если $x = x'$ и $y \neq y'$, то из (\ref{eq:vw})
имеем $\rho(r,r') \geq \rho(x+y,x'+y')+\rho(y,y') = 2\rho(y,y') \geq 4$.
Случай $x = x'$, $z \neq z'$ аналогичен.
Таким образом, при разных выборах $x$, $y$ и $z$  мы будем получать разные слова
с попарным расстоянием не менее $4$ друг от друга. Число всех таких слов равно
$|C_8|\cdot|B_8|^2 = 2^{18}$.
\proofend

Поскольку $F$ и $D$ есть линейные подпространства и, очевидно, $F \subset D$,
мы можем разбить $D$ на $64$ смежных класса по $F$, обозначим их $F_1,F_2,\ldots,F_{64}$.
Из леммы~\ref{l:pho4} следует, что окрестности этих смежных классов не пересекаются и,
воспользовавшись леммой~\ref{l:union}, мы можем строить совершенные раскраски
графа $\HHb{24}$ с параметрами вида
$((20+3i,256-3i)(3i,276-3i))$, $i=1,\ldots,64$.
Чтобы накрыть больший спектр параметров,
нам понадобится ещё один код:
$$L=\{(x,y,y+z)|x,y \in B, z \in C_1 \}$$
\begin{lemma}\label{l:8}
  Характеристическая функция $\chi_L$ есть совершенная раскраска графа $\HHa{24}$
  с параметрами $((28,248)(8,268))$.
\end{lemma}
\proofr
Представим $L$ в виде
$$L=\{(x,w)|x \in B, w \in C_{16} \}$$
где
$C_{16}=\{(y,y+z)|y \in B, z \in C_1 \}$ есть $(16,2^{11},4)$-код.

1) Рассмотрим кодовую вершину $(x,w)$ из $C_{16}$.
Слова кода $C_{16}$,
смежные вершине $(x,w)$ в $\HHa{24}$, имеют вид $(x+e,w)$,
где $e$ -- произвольное слово с двумя единицами. Поскольку таких слов $e$ ровно $28$,
каждая кодовая вершина смежна с $28$ кодовыми и, следовательно, с $248$ некодовыми.

2) Рассмотрим некодовую вершину $(x,w)$ графа $\HHa{24}$.
Если $x \not\in B$, то кодовые вершины графа $\HHa{24}$,
смежные вершине $(x,w)$, имеют вид $(x+e,w+e')$, где слова $e$ и $e'$
содержат по одной единице ровно. Слово $e$ можно выбрать восемью способами,
а $e'$ -- не более чем одним
(иначе в $C_{16}$ найдутся два слова на расстоянии $2$ друг от друга).
Таким образом, число кодовых вершин, смежных $(x,w)$, не превосходит $8$.
Если $x \in B$, то $w \not\in C_{16}$, и кодовые вершины,
смежные вершине $(x,w)$, имеют вид $(x,w+e'')$, где $e''$ содержит ровно две единицы.
Число способов выбора $e''$ -- не более $8$
(иначе в $C_{16}$ найдутся два слова на расстоянии $2$ друг от друга).
Таким образом, каждая некодовая вершина смежна не более чем с $8$ кодовыми.
С другой стороны, как следует из п.1, число рёбер,
соединяющих кодовые и некодовые вершины,
равно $2^{18}\cdot 248$ (где $2^{18}$ -- число всех кодовых вершин), что совпадает
с $(2^{23}-2^{18})\cdot 8$,
где  $2^{23}-2^{18}$ -- число всех некодовых вершин графа $\HHa{24}$.
Получается, что каждая некодовая вершина смежна ровно с $8$ кодовыми и,
следовательно, с $268$ некодовыми.
\proofend
Рассмотрим множество
$$N=\{(x+00000001,y+00000001,z+00000001)|x,y,z \in B)\subset\HHb{24}$$
и разобьём его на $8$ смежных классов $L_1$, \ldots, $L_8$ по $L$.
Поскольку расстояние от $D$ до $N$ равно $3$, все множества
 $\Omega(F_1)$, $\Omega(F_2)$, \ldots, $\Omega(F_{64})$, $L_1$, \ldots, $L_8$
 попарно не пересекаются и, применяя лемму~\ref{l:union}, получаем следующее:

 \begin{lemma}
   Для любых $i\in\{0,1, \ldots, 64\}$, $j\in\{0,1, \ldots, 8\}$, $0<i+j<72$,
  характеристическая функция объединения
  $i$ множеств из
  $\Omega(F_1)$, $\Omega(F_2)$, \ldots, $\Omega(F_{64})$
  и $j$ множеств из $L_1$, \ldots, $L_8$ есть совершенная раскраска с параметрами
  $((20+3i+8j, 256-3i-8j)( 3i+8j, 276-3i-8j))$.
 \end{lemma}
 Поскольку в виде $3i+8j$ могут быть представлены все числа от $3$ до $128$,
 кроме $4$, $5$, $7$, $10$, $13$,  теорема~\ref{th:main} доказана.

 \section[Proof of Theorem~2. The nonexistence]{Доказательство теоремы~\ref{th:main2}. Несуществование}
 \label{s:nonexi}
В этом разделе мы докажем
несуществование совершенных раскрасок с параметрами
$((20+c,256-c)(c,276-c))$ при $c=1,2,4,5,7$.
Множество $S \subset V(\HHa{24})$ назовём \emph{сферой}, если оно состоит из всех $24$ вершин
на расстоянии один от некоторой вершины (\emph{центра} сферы) из $V(\HHb{24})$.

В доказательстве следующей леммы будем обозначать двоичные слова длины $24$ в виде $V_{24}$,
где $V$ -- множество ненулевых позиций слова.

\begin{lemma} \label{l:12345} Предположим, что характеристическая функция $\chi_C$
множества $C \subset V(\HHa{24})$
является совершенной раскраской графа $\HHa{24}$ с параметрами $((20+c,256-c)(c,276-c))$.
Если $c\leq 7$, то $C$ есть объединение сфер.
\end{lemma}
\proofr

Мы рассмотрим только два случая $((25,251)(5,271))$ и $((27,249)(7,269))$,
поскольку остальные доказываются аналогично.

В первом случае $c=5$.
Возьмём произвольную вершину $v$ из $C$ и покажем, что она принадлежит
некоторой сфере, полностью содержащейся в $C$. Без потери общности будем считать, что
$v=000000000000000000000000$.

Пару координат $\{i,j\}$ от $1$ до $24$ назовём
\emph{кодовой}, если $\{i,j\}_{24}\in C$.
Из параметров совершенной раскраски следует, что кодовых пар ровно $25$.
Кроме того, поскольку $c=5$,

(*) некодовая пара может пересекаться не более чем с четырьмя кодовыми
(пятым соседом является вершина $v$).

Рассмотрим случаи:

1) Если в некоторой $i$-й координате пересекаются $23$ кодовые пары, то соответствующие
слова вместе с $v$ образуют сферу, что доказывает утверждение.

2) Если в некоторой $i$-й координате пересекаются от $5$ до $23$ кодовых пар,
то $i$ содержится в некоторой некодовой паре, что противоречит (*).

3) Если в некоторой $i$-й координате пересекаются ровно $4$ кодовые пары, то найдётся
кодовая пара $\{j,k\}$, не пересекающаяся ни с одной из них.
Тогда некодовая пара $\{i,j\}$ противоречит (*).

4) Пусть ни в одной координате не пересекаются $4$ или более кодовых пар.
Поскольку кодовых пар больше чем $24$, найдётся координата $i$, в которой
пересекаются $3$ кодовые пары
$\{i,j_1\}$, $\{i,j_2\}$, $\{i,j_3\}$. Легко посчитать,
что среди оставшихся $20$ координат найдётся такая $j$,
в которой пересекаются минимум две кодовые пары.
Тогда пара $\{i,j\}$ противоречит (*).

Для случая $c=5$ утверждение леммы доказано.

Рассмотрим случай $((27,249)(7,269))$, т.\,е. $c=7$.
В целом смысл доказательства похож на предыдущий рассмотренный случай,
только теперь мы рассматриваем не одну, а две соседних вершины первого цвета
и показываем, что при любом способе раскраски окрестности этих двух вершин
выполнение условия $a=27$ влечёт противоречие с условием $c=7$
(предположив от противного, что $C$ не есть объединение сфер).

Через $C'$ обозначим объединение сфер, полностью содержащихся в $C$;
а через $C''$ -- его дополнение до $C$.
Нам нужно доказать, что $C''$ -- пустое множество.
Предположим противное.
Обозначим через $k$ максимальную мощность пересечения кода $C$ и сферы,
содержащей хотя бы одну вершину из $C''$.
Поскольку такая сфера обязана содержать вершину не из $C$, имеем
\begin{equation}\label{eq:k7}
k\leq c=7.
\end{equation}
Без потери общности можем считать, что
 ${\bar 0}=000000000000000000000000\in C''$ и $\{1,2\}_{24},\{1,3\}_{24},\ldots,\{1,k\}_{24}\in C$.

Из слов кода $C$, принадлежащих окрестности ${\bar 0}$ или $\{1,2\}_{24}$, как из строк составим три
матрицы. Строки матрицы $A_1$ есть векторы веса $2$, отличные от $\{1,2\}_{24}$ и не лежащие
в окрестности $\{1,2\}_{24}$. Векторы веса $2$ ($4$), принадлежащие окрестности $\{1,2\}_{24}$,
составляют матрицу $A_2$ (соответственно $A_3$). Сами векторы ${\bar 0}$ и $\{1,2\}_{24}$ не включены
ни в одну матрицу. Матрицы $A_1$ и $A_2$ (как и $A_2$, $A_3$) состоят в совокупности из
$a-1=26$
строк, которые вместе с $\{1,2\}_{24}$ (соответственно с ${\bar 0}$) составляют пересечение
$C$ и окрестности ${\bar 0}$ (соответственно $\{1,2\}_{24}$).
Обозначив высоту матрицы $A_i$ через $h_i$, имеем:
$$h_1+h_2= h_2 + h_3 = 26.$$
Кроме того, из определения $k$ имеем $$h_2\leq 2k-4,$$ поскольку строки $A_2$
имеют вид $\{1,j\}_{24}$ (не более $k-2$, поскольку вместе с ${\bar 0}$ и $\{1,2\}_{24}$
принадлежат одной сфере) или $\{2,j\}_{24}$ (аналогично).

Так как строки $A_1$ и $A_2$ содержат по две единицы,
число единиц в этих двух матрицах равно  $52$.
Значит, в каком-то из столбцов содержится (суммарно в двух матрицах) не менее трёх единиц,
что по определению $k$ означает
\begin{equation}\label{eq:k4}
k\geq 4.
\end{equation}
Строки $A_3$ содержат по четыре единицы, значит
в совокупности в трёх матрицах имеем
\begin{equation}\label{eq:111}
2h_1+2h_2+4h_3 = 2\cdot 26 + 2(26-h_2)+2h_3 \geq 104 - 2(2k-4)+2h_3 = 112-4k+2h_3
\end{equation}
единиц ($2h_3$ единиц в первых двух столбцах $A_3$ выделяем отдельно, чтобы потом сократить).

Оценим это число с другой стороны. Назовём столбцы с номерами не больше (больше)
$k$ левой (соответственно правой) частью матрицы.

(a) \emph{Левая часть матрицы $A_1$ содержит не более $(k-2)(k-1)$ единиц.} Действительно,
первые два столбца этой матрицы нулевые, а каждый из оставшихся содержит не более
$k-1$ единиц, по определению $k$.

(b) \emph{Левая часть матрицы $A_2$ содержит не более $4(k-2)$ единиц.} Это следует
из уже рассмотренного состава матрицы.

(c) \emph{Правые части всех трёх матриц содержат в совокупности не более $(24-k)(7-k)$ единиц.}
Действительно, если в $j$-м столбце, $j>k$, матрицы $A_1$, $A_2$, $A_3$ содержат больше чем
$(7-k)$ единиц, то слово $\{1,j\}_{24}\not\in C$ имеет более $7$ соседей из $C$ (соответствующие строки
плюс ${\bar 0}$, $\{1,2\}_{24}$, $\{1,3\}_{24}$, \ldots, $\{1,k\}_{24}$), что противоречит параметру $c=7$.

Далее рассмотрим отдельно два случая:

I. $\{1,2\}_{24}\in C''$. В этом случае имеем также следующее:

(d) \emph{Левая часть матрицы $A_3$ содержит не более $2h_3+(k-2)(k-1)$ единиц.}
(Аналогично (a), только первые два столбца состоят из единиц.)

Итого, оценив суммарное число единиц в $A_1$, $A_2$, $A_3$ согласно (a)-(d)
и учитывая $(\ref{eq:111})$, имеем
$$112-4k+2h_3 \leq (k-2)(k-1)+4(k-2)+2h_3+(k-2)(k-1)+(24-k)(7-k).$$
То есть,
$$ 3k^2-29k+52 \geq 0,$$
что неверно при значениях $k=4,5,6,7$, удовлетворяющих (\ref{eq:k7}) и (\ref{eq:k4}).
Получили противоречие.

II. $\{1,2\}_{24}\in C'$. Тогда $C$ включает в себя сферу с центром $\{1,2,j\}_{24}$ для
некоторого $j$. Следовательно, матрица $A_3$ содержит все строки вида $\{1,2,j,i\}_{24}$,
$i\in \{3,\ldots,24\}\setminus\{j\}$, и, как следствие, каждый столбец $A_3$ содержит
не менее одной единицы. Учитывая (c), получаем

(e) \emph{Правые части $A_1$ и $A_2$ содержат в совокупности не более $(24-k)(6-k)$ единиц.}

Как мы уже отмечали, число единиц в матрицах $A_1$ и $A_2$ равно $52$.
С другой стороны, как следует из
(a), (b) и (e), оно не превосходит
$$ (k-2)(k-1) + 4(k-2) + (24-k)(6-k) = 2k^2 - 29k +138. $$
То есть получаем $2k^2 - 29k + 86 \geq 0$, что не верно для $5\leq k \leq 7$,
но справедливо при $k=4$.

Этот последний оставшийся подслучай рассмотрим отдельно.
Имеем: $\{1,2\}_{24}$, $\{1,3\}_{24}$, $\{1,4\}_{24}$ $\in C$, причём
$\{1,2\}_{24}$ принадлежит некоторой сфере, полностью лежащей в $C$.
Центр этой сферы имеет вид $\{1,2,j\}_{24}$ (он не может содержать одну единицу, поскольку
$\bar 0 \in C''$). Поскольку тогда $\{1,j\}_{24}\in C$, $j$ есть $3$ или $4$.
Рассмотрим случай $j=3$ (второй аналогичен). Покажем, что

(**) $\{1,4\}_{24}\in C''$. Пусть, от противного,
вершина $\{1,4\}_{24}$ принадлежит некоторой сфере, полностью лежащей в $C$.
Аналогично уже рассмотренному случаю с вершиной $\{1,2\}_{24}$, центр этой сферы
обязан иметь вид $\{1,j,4\}_{24}$, где $j$ есть $2$ или $3$. Но тогда сферы
с центрами $\{1,2,3\}_{24}$ и $\{1,j,4\}_{24}$ пересекаются. Точка пересечения
принадлежит $C$ и имеет не менее $45$ соседей из $C$ ($1+45$ есть мощность объединения
двух сфер), что противоречит параметру $a=27$. Утверждение (**) доказано.

Теперь, поменяв местами вторую и четвёртую координаты,
мы приходим к уже рассмотренному случаю I.
\proofend

\begin{lemma}\label{l:::3}
Пусть $C \subset V(\HHa{24})$ есть объединение сфер, и характеристическая функция $\chi_C$
есть совершенная раскраска графа $\HHa{24}$ с параметрами $((20+c,256-c)(c,276-c))$.
Тогда либо $c$ кратно $3$, либо $c \geq 25$.
\end{lemma}
\proofr
Если $C$ включает в себя две пересекающиеся сферы, то их объединение содержит $46$ вершин,
а любой (из двух) общий элемент смежен с остальными $45$ вершинами объединения. Отсюда
$20+c \geq 45$, т.\,е. $c \geq 25$.

Если, в противном случае, $C$ состоит из непересекающихся сфер,
то каждая вершина графа $\HHa{24}$
не из $C$ смежна с тремя вершинами из каждой соседней сферы, откуда $c\equiv 0 \bmod 3$.
\proofend
Поскольку значения $c=1,2,4,5,7$ противоречат леммам~\ref{l:12345} и~\ref{l:::3},
теорема~\ref{th:main2} доказана.
Кроме того, можно сделать вывод, что совершенная раскраска с параметрами $((23,253)(3,273))$
единственна с точностью до автоморфизма графа.
\section[Conclusion]{Заключение}

В заключение приведём таблицу всевозможных значений параметра $c$ от $1$ до $128$.
Знак "$-$"\ означает несуществование совершенной раскраски с параметрами
$(20+c,256-c;c,276-c)$ в $\HH{24}$,
"$+$"\ --
существование, "$?$"\  -- вопрос существования открыт.
Если плюс обведён в кружок, то раскраску с данными
параметрами можно построить объединением смежных классов окрестности кода Голея,
если в прямоугольник --
объединением смежных классов линейного кода $L$.

$$\scriptstyle
\def\cpl{\makebox[0mm][c]{$\bigcirc$}\makebox[0mm][c]{$+$}}
\def\qpl{\framebox[4.5mm]{$+$}}
\def\qcpl{\framebox[4.5mm]{\makebox[0mm][c]{$\bigcirc$}\makebox[0mm][c]{$+$}}}
\begin{array}{|@{\,}c@{\,}|@{\,}c@{\,}|@{\,}c@{\,}|@{\,}c@{\,}|@{\,}c@{\,}|@{\,}c@{\,}|@{\,}c@{\,}|@{\,}c@{\,}|@{\,}c@{\,}|@{\,}c@{\,}|@{\,}c@{\,}|@{\,}c@{\,}|@{\,}c@{\,}|@{\,}c@{\,}|@{\,}c@{\,}|@{\,}c@{\,}|@{\,}c@{\,}|@{\,}c@{\,}|@{\,}c@{\,}|@{\,}c@{\,}|@{\,}c@{\,}|@{\,}c@{\,}|@{\,}c@{\,}|@{\,}c@{\,}|@{\,}c@{\,}|@{\,}c@{\,}|@{\,}c@{\,}|@{\,}c@{\,}|@{\,}c@{\,}|@{\,}c@{\,}|@{\,}c@{\,}|@{\,}c@{\,}|}
\hline
1 & 2 & 3 & 4 & 5 & 6 & 7 & 8 & 9 & 10 & 11 & 12 & 13 & 14 & 15 & 16 \\
-&-&\cpl&-&-&\cpl&-&\qpl&\cpl&?&+&\cpl&?&+&\cpl&\qpl\\ \hline
 17 & 18 & 19 & 20 & 21 & 22 & 23 & 24 & 25 & 26 & 27 & 28 & 29 & 30 & 31 & 32 \\
+&\cpl&+&+&\cpl&+&+&\qcpl&+&+&\cpl&+&+&\cpl&+&\qpl\\ \hline
33 & 34 & 35 & 36 & 37 & 38 & 39 & 40 & 41 & 42 & 43 & 44 & 45 & 46 & 47 & 48 \\
\cpl&+&+&\cpl&+&+&\cpl&\qpl&+&\cpl&+&+&\cpl&+&+&\qcpl\\ \hline
 49 & 50 & 51 & 52 & 53 & 54 & 55 & 56 & 57 & 58 & 59 & 60 & 61 & 62 & 63 & 64 \\
+&+&\cpl&+&+&\cpl&+&\qpl&\cpl&+&+&\cpl&+&+&\cpl&\qcpl\\ \hline
65  &   66  &   67  &   68  &   69  &   70  &   71  &   72  &   73  &   74  &   75  &   76  &   77  &   78  &   79  &   80  \\
+&\cpl&\cpl&+&\cpl&\cpl&+&\qcpl&\cpl&+&\cpl&\cpl&+&\cpl&\cpl&\qpl\\ \hline
81  &   82  &   83  &   84  &   85  &   86  &   87  &   88  &   89  &   90  &   91  &   92  &   93  &   94  &   95  &   96  \\
\cpl&\cpl&+&\cpl&\cpl&+&\cpl&\qcpl&+&\cpl&\cpl&+&\cpl&\cpl&+&\qcpl\\ \hline
97  &   98  &   99  &   100 &   101 &   102 &   103 &   104 &   105 &   106 &   107 &   108 &   109 &   110 &   111 &   112 \\
\cpl&+&\cpl&\cpl&+&\cpl&\cpl&\qpl&\cpl&\cpl&+&\cpl&\cpl&+&\cpl&\qcpl\\ \hline
113 &   114 &   115 &   116 &   117 &   118 &   119 &   120 &   121 &   122 &   123 &   124 &   125 &   126 &   127 &   128 \\
+&\cpl&\cpl&+&\cpl&\cpl&+&\qcpl&\cpl&+&\cpl&\cpl&+&\cpl&\cpl&\qpl\\ \hline
\end{array}
$$

\providecommand\href[2]{#2} \providecommand\url[1]{\href{#1}{#1}}

\parbox[t]{30em}{
Адрес автора:\\[1ex]
Институт математики\\
им. С. Л. Соболева СО РАН,\\
пр-т Академика Коптюга, 4,\\
630090 Новосибирск,\\
Россия.\\[0.5ex]
E-mail: \texttt{krotov@math.nsc.ru}
}
\hfill
\parbox[t]{20em}{Статья поступила\\ 19 марта 2008\,г.}